\documentclass[12pt]{article}
\usepackage{amsmath,amsfonts,amssymb,amsthm}
\usepackage{graphicx}
\usepackage{caption}
\usepackage{subcaption}
\usepackage{mathtools}
\usepackage{cite}
\usepackage{hyperref}
\usepackage{enumerate}

\makeatletter
\@addtoreset{equation}{section}
\makeatother

\def\eqref#1{(\ref{#1})}

\newcommand{\per}{\mathbb{P}\mathrm{er}}
\newcommand{\length}{\mathrm{length}}
\newcommand{\norm}[1]{\left\lVert#1\right\rVert}
\let\phi=\varphi 
\let\epsilon=\varepsilon 
\newcommand{\comp}{\texttt{Comp}}
\newcommand{\diff}{\texttt{Diff}}
\newcommand{\teich}{\texttt{Teich}}
\newcommand{\C}{\mathcal{C}}
\newcommand{\Lie}{\operatorname{Lie}}
\newcommand{\Hom}{\operatorname{Hom}}
\newcommand{\End}{\operatorname{End}}

\newcommand{\Tw}{\operatorname{Tw}}

\title{Sub-twistor metrics}
\author{Dmitrii Korshunov}

\theoremstyle{definition}
\newtheorem{defn}{Definition}[subsection]
\newtheorem{example}{Example}[subsection]
\newtheorem{lemma}{Lemma}[subsection]
\newtheorem{proposition}{Proposition}[subsection]
\newtheorem{theorem}{Theorem}[section]
\newtheorem{remark}{Remark}[subsection]

\begin{document}
\maketitle

\begin{abstract}
We consider a natural distance function on the period space of a hyperk\"ahler manifold associated to non-holonomic constraints imposed by twistor lines. These metrics were introduced by Verbitsky in the context of the global Torelli theorem for hyperk\"aler manifolds. We show that they are Finsler and explicitly describe  their Finsler norms. To achieve this goal we consider a class of distance functions (called here sub-conic metrics) that slightly generalize sub-Riemann metrics. The main technical result is the statement that every sub-conic metric is sub-Finsler. The content and methods of the paper lie within basic metric geometry. The complex geometrical background and motivation are isolated in the appendix.
\end{abstract}
\tableofcontents

\section{Introduction}
Given a manifold $M$ endowed with a riemannian form $g\in \Gamma(M,S^2T^*M)$ the length of a piecewise smooth path $\gamma: [0,t]\to M$ is classically defined by

\begin{equation}
    \label{distance}
\length_g(\gamma):=\int_0^t\sqrt{g(\dot\gamma(s),\dot\gamma(s))}ds
\end{equation}

where $\dot\gamma$ is a velocity vector, that is $\gamma_*(\frac{d}{ds})$, which is well defined almost everywhere.

The Riemannian distance between two points $x$ and $y$ in $M$ is then introduced as

\begin{equation}\label{eqn:distance}
d_g(x,y):=\lim_{\gamma\in H}\length(\gamma)
\end{equation}

where infimum is taken over the set $H$ of all piecewise smooth paths $\gamma$ that connect $x$ and $y$.

There are two straightforward ways to generalize this notion of distance. We can modify the way we measure lengths of paths by substituting (imposing suitable regularity conditions) a quadratic norm $\sqrt{g(\cdot,\cdot)}$ by an arbitrary norm $\norm\cdot$. The metrics defined in such way are called Finsler. The second possible direction is to restrict the class $H$ of paths over which the limit is taken. For example, one can consider only those piecewise smooth paths whose velocity vectors $\dot\gamma$, whenever defined, lie in a fixed sub-bundle $D$ of a tangent bundle of $M$.  There is a huge body of literature devoted to the topic. Being approached from different directions it is known under various names: non-holonomic dynamics \cite{vershik}, Carnot-{C}arath\'{e}odory geometry \cite{gromov}, sub-Riemannian geometry \cite{bellaiche} and geometric control theory \cite{Agrachev, jean}. 
The simultaneous generalization, when we restrict the class of paths and define a length functional using an arbitrary norm, is called sub-Finsler metric. The interest in sub-Finsler metrics is relatively recent \cite{berestovski,donna,Bre, subfinsler}.

In this paper we will be concerned with the generalization of the second type. Our object of interest will be the period space $\per$ of a hyperk\"ahler manifold \cite{Ver, Huy}. There is a rich set of distinguished two-dimensional submanifolds of $\per$ called {\it twistor spheres} and we define the sub-twistor metric by the formula \eqref{distance} taking the limit over the set of piecewise smooth paths tangent to twistor spheres almost everywhere. We will be referring to such paths as {\it twistor paths}. There are enough twistor spheres in $\per$ for this metric to be well defined, namely any two points can be connected by a twistor path.

All tangent planes of twistor spheres passing through a point  sweep out a cone of possible velocities of  twistor paths. Under a mild regularity assumption, we show that any metric defined by such ``distribution of cones'' is sub-Finsler (in a sense described below). The sub-twistor metric  turns out to be a genuine Finsler metric. In particular, it induces the underlying topology of the manifold $\per$.

The sub-twistor metric was introduced by Verbitsky \cite{Ver} in course of his proof of the global Torelli theorem for hyperk\"ahler manifolds. Loosely speaking, the theorem states that a certain map (the period map) between two manifolds --- a connected component of the birational Teichm\"uller space and the period space $\per$ --- is a homeomorphism. The period map is a local isometry with respect to the sub-twistor metric by the lifting property of twistor spheres. By the path-lifting property for the period map and simply-connectedness of $\per$ it follows that the period map is a homeomorphism.

The result of the paper and methods are those of metric geometry. All the algebro-geometric motivation is isolated in the appendix. 

\section{Sub-conic metrics}
\subsection{Sub-conic and sub-Finsler metrics}
\begin{defn}\label{conedef}
A subset $D\subset TM$  of the total space of the tangent bundle of a manifold $M$  is called a {\it distribution of cones} if it is a union of graphs of smooth vector fields $X$, and is closed with respect to multiplication by smooth functions $f\in C^\infty(M)$ in the following sense: if the graph of $X$ is a subset of $D$ then the graph of $fX$ is also a subset of $D$.
\end{defn}

By a slight abuse of notation we will write $X\in D$ when the graph of a vector field $X$ is a subset of $D$.
By {\it the value of a vector field} $X$ at a point $p$ we mean the vector $X_p\in T_p M$ that is a value of $X$ considered as a section of $TM\to M$.

\begin{remark}
The set of values of all vector fields $X\in D$ at a point $p$ is a cone in $T_pM$, i.e. a set of vectors invariant with respect to multiplication by scalars. We  denote this set by $D_p\subset{T_pM}$. 
\end{remark}

\begin{defn}
A distribution $S$ on a manifold $M$ is a distribution of cones with the following additional property: if the graphs of vector fields $X$ and $Y$ are subsets of $S$ then the graph of the vector field $fX+gY$, where $f,g\in C^\infty(M)$ is also a subset of $S$. We will write $X\in S$ to indicate that the graph of $X$ is a subset of $S$.
\end{defn}
\begin{remark}
At any point $p$ the set of all values of vector fields from a distribution $S$ is the linear subspace which we will denote $S_p\subset T_pM$. Notice that we do not require a distribution to be a subbundle, thus $\dim S_p$ need not be constant. A distribution is just a distribution of cones that is closed with respect to addition of vector fields.
\end{remark}

\begin{defn}
A set of vector fields $S$ is said to satisfy {\it the Chow condition} if the Lie algebra generated by sections $X\in S$ is all $\Gamma(M,TM)$. In other words, nested Lie brackets of sections from $S$ generate $\Gamma(M,TM)$. That is, at any point $p$ linear combinations of values of $X_i$, $[X_1,X_2]$, $[[X_1,X_2],X_3]$,  $[[X_1,X_2],[X_3,X_4]]$ etc., $X_i\in D$, span $T_pM$.
\end{defn}

Given a distribution of cones, we can take a fiberwise linear span of cones in each $T_pM$, which we will refer to as the {\it distribution associated to a cone}.

\begin{defn}
A distribution $\langle D\rangle$ associated to a distribution of cones $D$ is the linear span of $D$, i.e. the set of all vector fields of the form $a_1X_1+a_2X_2+\dots+a_mX_m$ where $a_i\in C^\infty(M)$ and $X_i\in D$.
\end{defn}

\begin{remark}
At a point $p$ the set of values  of vector fields from $\langle D\rangle$ coincides with the linear span of $D_p$, that is, $\langle D \rangle_p=\langle D_p \rangle$. Indeed, let $X_p^i\in D_p$ and $w=\sum_i \lambda_i X_p^i$. Each $X_p^i$ is a value of some vector field $X^i$ from $D$, and $W=\sum_i \lambda_i X^i$ is still a vector field from $D$ whose value at $p$ is $w$. On the other hand, the value of $\sum_i f_i X^i$ at $p$ is  $\sum_i f_i(p) X^i_p\in \langle D_p\rangle$.
\end{remark}

To any distribution of cones $D$ (in particular, any distribution) we can associate a distance function by restricting the class of paths to those with velocity vectors lying in $D$.

\begin{defn}
A {\it $D$-admissible path}  is a piecewise smooth map $[0,t]\to M$ such that at any point $p=\gamma(s)$ where velocity vector $\dot\gamma:=\gamma_*(\frac{d}{dt})$ is defined, $\dot\gamma(s)\in D_p$.
\end{defn}

\begin{defn}
The {\it sub-conical metric}, associated to a triple $(M,g,D)$, where $g$ is a Riemannian metric on $M$ and $D$ a distribution of cones, is a distance function on $M$ defined by

$$d_D(x,y):=\inf_\gamma \int_0^t \sqrt{g(\dot\gamma(s),\dot\gamma(s))}ds$$

where infimum is taken over the set of all $D$-admissible paths $\gamma$  such that $\gamma(0)=x$ and $\gamma(t)=y$.
\end{defn}

\begin{figure}
\centering
\begin{minipage}{.5\textwidth}
  \centering
  \includegraphics[width=.82\linewidth]{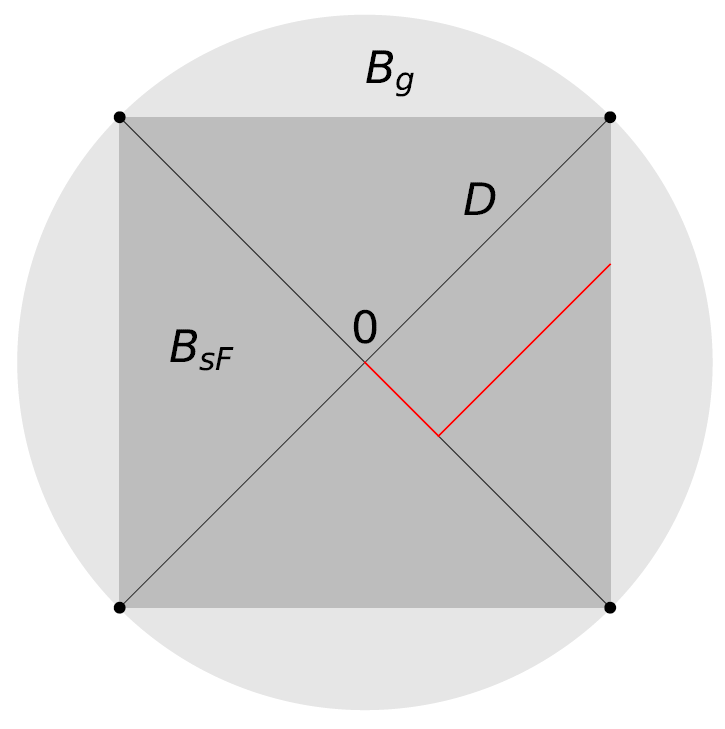}
  \captionof{figure}{Finsler norm of Manhattan metric}\label{manhattan}
  \label{fig:convex_hull}
\end{minipage}%
\begin{minipage}{.5\textwidth}
  \centering
  \includegraphics[width=.9\linewidth]{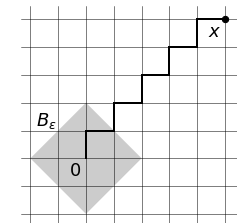}
  \captionof{figure}{Manhattan metric}
  \label{fig:taxicab.png}
\end{minipage}
\end{figure}

\begin{example}
Consider $\mathbb{R}^2$ equipped with the standard euclidean metric. Let $X=\frac{\partial}{\partial x}$ and $Y=\frac{\partial}{\partial y}$ be the standard coordinate vector fields. Define a distribution of cones by
$D=\big\{fX,gY: f,g\in C^\infty(\mathbb{R}^2) \big\}$
. In this case admissible paths are those that are piece-wise horizontal or vertical. The sub-conical distance associated to this data is the well-known Manhattan (or $\ell_1$) distance. It is also a Finsler distance on $\mathbb{R}^2$, with the norm given by $\nu(a X+ b Y)=|a|+|b|$. The unit ball is shown in Figure \ref{manhattan}.
\end{example}

This turns out to be a general phenomenon: any sub-conical metric associated to a conical distribution $D$ is equal to a sub-Finsler metric, defined on the fiber-wise linear span of $D$.

A sub-Finsler structure is a distribution on a manifold equipped with a norm in each fiber. However, different regularity assumptions are possible. In the Finsler geometry literature it is often imposed that Finsler metric is smooth away from the zero section. This definition is too restrictive for our purposes and we require a Finsler norm to be a merely continuous function (in accordance, e.g. with \cite{berestovski}).
The interest in sub-Finsler metrics is quite recent in comparison to Finsler and sub-Riemannian geometries and stems mainly from the problems of geometric group theory \cite{Bre} and geometric  control theory \cite{subfinsler}. The importance of sub-Finsler metrics is highlighted by a theorem of Berestovski \cite{berestovski}, which says that any left-invariant inner metric on a Lie group is sub-Finsler.
\begin{defn}
A sub-Finsler structure on a manifold $M$ equipped with a distribution $S$ viewed as a subset of the total space of $TM$ is a continuous function $\nu$ on $S$, such that the restriction to any fiber $S_p$ is a norm.
\end{defn}

Now we can define the sub-Finsler distance function:

\begin{defn}
Given a sub-Finsler structure $(M,S,\nu)$, the corresponding sub-Finsler distance is defined by 

$$d_{sF}(x,y):=\inf_\gamma\int_0^t\nu(\dot\gamma(s))ds$$

where infimum is taken over all $S$-admissible paths, connecting $x$ and $y$.
\end{defn}

\begin{defn}\label{assfinsler}
Let $M$ be a manifold equipped with a distribution of cones $D$ and a Riemannian metric $g$. Define a norm on  $\langle D \rangle$ by the formula

$$\norm v_{g,D}=\inf\sum_i\lambda_i$$

where infimum is taken over all possible expressions of $v$ as a convex combination of vectors from $D_p$

$$v=\sum \lambda_iX_i$$

where $\lambda_i>0$, $X_i\in D_p$ and $g(X_i,X_i)=1$.

We will be referring to this norm as {\it the sub-Finsler norm associated to a distribution of cones $D$ and a Riemannian metric $g$}.
\end{defn}

\begin{remark}
$\norm \cdot_{g,D}$ is indeed a norm, that is 
\begin{enumerate}
  \item $\norm{\lambda v}_{g,D}= |\lambda|\cdot\norm v_{g,D}$
  \item $\norm v_{g,D}=0$ if and only if $v=0$
   \item $\norm{v+w}_{g,D}\le \norm v_{g,D} +  \norm v_{g,D} $
\end{enumerate}
\end{remark}

\begin{remark}\label{conedir}
If $v\in D_p$ then $\norm v_{g,C}=\sqrt{g(v,v)}$. In other words, if $v\in D_p$ then its sub-Finsler norm is equal to its Riemannian norm. For $v\not\in D_p$, however, $\norm{v}_{g,C}$ is generally less than $\sqrt{g(v,v)}$.
\end{remark}

\begin{remark}
\label{norm_closed}
Any norm on $\mathbb{R}^n$ is defined by its closed unit ball, a centrally symmetric bounded convex body with nonempty interior. Any bounded closed centrally symmetric convex body $B$ that is not a subset of a proper linear subspace is a unit ball of the norm $\norm v_B:=\inf\{\lambda>0: \lambda v\in B\}$. With this correspondence in mind, the unit ball for the norm $\norm{\cdot}_{g,D}$ can be described as
$$B_{sF}=\overline{\mathrm{Conv}}(D_p\cap B_g)$$
where $\overline{\mathrm{Conv}}$ denotes the closure of the convex span of a subset. Indeed, the unit ball $B_{sF}$ is the closure of the set of all convex combinations (i.e. linear combinations with coefficients positive real numbers summing up to $1$) of vectors from $D_p$ of Riemannian norm $\le 1$. Figure \ref{manhattan} illustrates the case of Manhattan metric.
\end{remark}

\begin{proposition}\label{cont_norm}
The norm $\norm{\cdot}_{g,D}$  is an upper semi-continuous function on $\langle D\rangle$ (viewed as a subset of the total space of $TM$).
\end{proposition}
\begin{proof}

First we prove that if an infimum of $\sum_i \lambda_i$ is $\epsilon$-approximately realized for a fixed finite set of tangent vectors at some point $p$ then it is $\epsilon$-approximately realized on $\dim\langle D_p \rangle$ linearly independent vectors from this set. Here we say that the norm is $\epsilon$-approximately realized when we chose a set of unit vectors $X_i$ and positive real numbers $\lambda_i$, $v=\sum_i\lambda_iX_i$, such that $|\norm{v}_{g,D}-\sum_i\lambda_i|<\epsilon$. Indeed, let $\{X_i\}_{i=1..m}$ be a set of unit vectors  from $D_p$ such that $v=\sum_{i=1}^m \lambda_iX_i$ for some $\lambda_1,\dots\lambda_m$. Consider the linear function $\mathbb R_+^m\to \mathbb R_+$ defined by $(\lambda_1,\dots \lambda_m)\to \lambda_1+\dots +\lambda_m$. The subspace of all $m$-tuples $(\lambda_1,\dots \lambda_m)$ that give a decomposition $v=\sum_{i=1}^m \lambda_iX_i$ is an affine subspace of $\mathbb R_+^m$. Unless this subspace is a point, the function $\lambda_1+\dots+\lambda_m$ takes its minimum at the boundary, that is when one of the $\lambda_i$ is zero. Thus we can remove corresponding vectors until the decomposition is unique, which implies that all the vectors are linearly independent.

Thus it is enough to evaluate the norm $\norm{\cdot}_{g,D}$ on linearly independent vector fields. Now, for any set $S$ of linearly independent unit vector fields $X_i\in D$  in a neighborhood of a point $p$, consider the function on $\langle D\rangle$:

$$f_S(v):=\sum_{X_i\in S}|\langle v,X_i\rangle|$$

This function is continuous. On the other hand, one now has

$$\norm v_{g,D}=\inf_S\big\{f_S(v)\big\}$$

Since the pointwise infimum of any set of continous functions is upper semi-continous, the $\norm v_{g,D}$ is upper semicontinous. 
\end{proof}

\begin{remark}
The norm as defined is not necessarily continuous. The simple example is the cone generated by $x\frac{\partial}{\partial x}$, $\frac{\partial}{\partial y}$ and $\frac{\partial}{\partial y}+\frac{\partial}{\partial x}$ on the standard $\mathbb{R}^2$. For the sequence of tangent vectors $v_n\in T_{-\frac{1}{n}}\mathbb R^2$ all equal to $\frac{\partial}{\partial x}$, one has $\norm {v_n}=1$, but $\norm{\frac{\partial}{\partial x}(0)}=1+\sqrt 2$.

However in all reasonable situations it is indeed  continuous. For example, when the distribution of cones $D$ is closed as a subset of $TM$.
\end{remark}

\begin{proposition}
Let $D$ be closed as a subset of   $TM$. Then the norm  $\norm{\cdot}_{g,D}$  is continuous.
\end{proposition}

\begin{proof}
In view of Proposition \ref{cont_norm}, it is left to prove that $\norm{\cdot}_{g,D}$ is lower semi-continuous. Let $v_i\in T_{p_i}M$ be a sequence of tangent vectors converging to some tangent vector $v\in T_pM$. Suppose that each $v_i\in TM$ is decomposed as $v_i=\sum_{j=1}^n \lambda_{ij} X_{ij}$ with $ \norm X_{ij}=1 $ realizing the norm. Now choose a compact neighborhood of $p$ and for each sequence $ X_{ij} $ pick a  subsequence that converges to some $ X_{j}\in T_pM$.  The decomposition $v=\sum_{j=1}^n \langle v, X_j \rangle X_j$  shows that the norm on $T_pM$ is not greater than the limit of norms on $T_{p_i}M$.
\end{proof}

\begin{remark} 
\label{fiberwise_clos}
Note that by Remark \ref{norm_closed} the metric only depends on the fiberwise closure of the distribution of cones. Thus from now on we will suppose, as part of the definition, that the fiberwise closure of all distributions of cones is closed in the above sense. This also implies continuity of the norm $\norm{\cdot}_{g,D}$ by the previous proposition. 
\end{remark}
    
\subsection{Sub-conical metrics are sub-Finsler}

The following theorem was proposed as a conjecture by S.V. Ivanov in a personal communication. It states that a sub-conical distance $d_D$, associated to a distribution of cones $(M,D,g)$, equals the sub-Finsler distance $d_{sF}$, associated to $(M,\langle D\rangle,\norm{\cdot}_{g,C})$.
In other words, we can always ``relax'' a distribution of cones to a distribution of linear subspaces, extending the class of admissible curves. 

\begin{theorem}
\label{mainthm}
Let $M$ be a Riemannian manifold with a metric $g$, $D$ a distribution of cones on $M$, such that $D$ satisfies the Chow condition. Denote by $d_D$ the sub-conical metric associated to the distribution of cones $(M,g,D)$ and let $d_{sF}$ be the sub-Finsler metric associated to the distribution $\langle D \rangle$ and the norm $\norm{\cdot}_{g,D}$ (Definition \ref{assfinsler}). Then

$$d_D=d_{sF}$$
\end{theorem}
\begin{proof}
Let $\gamma$ be a $D$-admissible path connecting points $x$ and $y$. Then by Remark \ref{conedir}, since $\dot\gamma\in D$ and $\norm{\dot\gamma}_{g,C}=g(\dot\gamma,\dot\gamma)$

$$\int_0^t\sqrt{ g(\dot\gamma,\dot\gamma)}ds=\int_0^t \norm{\dot\gamma}_{g,C}ds$$

Hence $d_D\ge d_{sF}$.

To  prove the opposite inequality $$d_D\le d_{sF}$$ it is enough to show that for any $\langle D\rangle$-admissible path $\gamma$ between $x$ and $y$ and any $\epsilon>0$ there exists a $D$-admissible path $\sigma$ connecting the same points, whose Riemannian length is $\epsilon$-close to the Finsler length of $\gamma$:

$$\Big\lvert\int_0^t \sqrt{(\dot{\sigma}(s),\dot{\sigma}(s))}ds -\int_0^t \norm{\dot\gamma(s)}ds\Big\lvert\le\epsilon$$
This is the statement of Proposition \ref{zig-zag}, proved later. Now, since $\epsilon$ can be chosen arbitrarily small

$$\inf_\sigma \int_0^t \sqrt{g(\dot{\sigma}(s),\dot{\sigma}(s))}ds \le \int_0^t \norm{\dot\gamma(s)}ds$$
and the inequality $d_D\le d_{sF}$ follows. The rest of this section will be devoted to the construction of the approximating $D$-admissible paths with the required property.
\end{proof}

We will achieve the approximation of a $\langle D\rangle$-admissible path by a $D$-admissible one in two steps. Firstly, we will use the Lie product formula for vector fields to get a uniform approximation of the flow of a vector field $X+Y$ by alternating concatenations of flows of $X$ and $Y$ (``zig-zag'') simultaneously estimating the lengths of zig-zags in terms of the Finsler length of $\gamma$. The second step will be the demonstration that zig-zag paths can be slightly modified to actually connect $x$ and $y$.

We use the notation $e^{tX}$ for the flow generated by  a vector field $X$. Thus $e^{tX}p:\mathbb{R}\to M$ is the orbit of a point $p$. It is well defined for small enough values of $t$.

\begin{subsection}{Zig-zag approximation}

In this subsection the notation $|x-y|$, where $x$ and $y$ are points on a manifold $M$, means the Riemannian distance between $x$ and $y$. When applied to real numbers it denotes the absolute value.
By $\length( \gamma)$ we will denote the Riemannian length of the the path $\gamma$. Thus $\length (e^{tX} p)$ is the Riemannian length of the orbit (up to the moment of time $t$) of a point $p$ with respect to the flow $e^{tX}$.

\begin{lemma}\label{abraham}
Let $X_1,\dots X_r$ be a finite set of smooth vector fields on a manifold $M$. Then for any $t$ for which the flows of $X_1,\dots X_r$ and  $X_1+\dots+X_r$ are defined and any $\epsilon>0$ there exists a  natural number $N>0$ such that
$$\Big\lvert \big(e^{\frac{t}{N}X_1}\dots e^{\frac{t}{N}X_r}\big)^Np-e^{t(X_1+\dots +X_r)}p\Big \rvert<\epsilon$$

\end{lemma}
\begin{proof}
This is Corollary 2.1.27 in Abraham-Marsden's ``Foundations of Mechanics" \cite{mechanics}.
\end{proof}

The following lemma bounds the rate at which lengths of integral curves of a vector field on a compact subset may differ if their initial points are close.

\begin{lemma}\label{stretch}
Let $X$ be a vector field on a compact neighborhood $U$ and $T>0$ a real number. Then there exists a constant $C$, depending on $ T $, $U$ and $X$, such that  for any $\epsilon>0$ and any $p,q\in U$ such that $|p-q|<\epsilon$, one has

$$\Big\lvert\length (e^{tX}p) - \length (e^{tX}q)\Big\rvert\le C\epsilon t$$
\end{lemma}

for all $t<T$ for which the  integral curve of the flow $e^{tX}$ remains in $U$.

\begin{proof}

Let $$C':= \sup_{p\in U,\,t\le T} \norm{De^{tX}(p)}$$
i.e. $C'$ is the supremum of operator norms of the differential of the flow at all points of  $U$ for every $0 \le t\le T$. Then $C'$ is finite by compactness of $U\times[0,T]$ and continuity of the operator norm. In other words, $C'$ is a local Lipschitz constant for  $e^{tX}$ for all $t\le T$. Since a locally Lipschitz function on a compact is Lipschitz, and a composition of locally Lipschitz functions is locally Lipschitz, one has for some constant $C$:
$$\Big\lvert  \norm{X(e^{tX}p)} - \norm{X(e^{tX}q)}  \Big\rvert\le C|p-q|\le C\epsilon$$
 Here $\norm{X}$ denotes the Riemannian norm of a vector $X$ (note that it is Lipschitz). Now by a standard inequality for the absolute value of an integral
 
$$\Big\lvert{\int_0^t \norm{X(e^{sX}p)}ds-\int_0^t \norm{X(e^{sX}q)}ds}\Big\rvert\le 
\int_0^t\Big\lvert  \norm{X(e^{sX}p)} - \norm{X(e^{sX}q)}  \Big\rvert ds\le$$

$$\le C\epsilon t$$

\end{proof}

\begin{lemma}\label{mainlemma}
For any smooth path $\gamma$ and vector fields $X_1,\dots X_r$ on a compact neighborhood $U$ of $p$, such that $\dot\gamma(t) = \sum X_i(\gamma(t))$ and any $t$ for which the flows of $X_1,\dots X_r$ and $\sum X_i$ are defined, the following holds:
$$\lim_{N\to\infty} \length((e^{\frac{t}{N}X_1}\dots e^{\frac{t}{N}X_r})^Np) = \int_0^t \sum_{i=1}^r\norm{X_i(\gamma(s))}ds$$

\end{lemma}

\begin{proof}
Pick $\epsilon>0$ and using Lemma \ref{abraham} choose $N$ large enough, such that 
$$\Big\rvert(e^{\frac{t}{N}X_1}\dots e^{\frac{t}{N}X_r})^Np-e^{t(X_1+ \dots + X_r)}p\Big\lvert<\epsilon$$

We will be referring to paths $Z^N(p,t):=(e^{\frac{t}{N}X_1}\dots e^{\frac{t}{N}X_r})^Np$ as zig-zags. Denote the points on $\gamma$, corresponding to the equidistant moments of time $$0, \frac{t}{N},\frac{2t}{N} \dots \frac{(N-1)t}{N}, t$$ by
$p_k:=\gamma\big(k\frac{t}{N}\big)$
and the points corresponding to the same moments of time on the approximating zig-zag curve

\[   \left\{
\begin{array}{ll}
      \widehat p_0=p_0=p \\
      \widehat p_{k+1}=e^{\frac{t}{N}X_1}\dots e^{\frac{t}{N}X_r}\widehat p_k\\
\end{array} 
\right . \]

Note that  $|p_k-\widehat p_k|\le \epsilon$ by the choice of $N$ since $p_k$ are points on $e^{t(X_1+ \dots + X_r)}p$ and $\widehat p_k$ are points on the zig-zag corresponding to the same moments of time $k\frac{t}{N}$. Indeed, Lemma \ref{abraham} asserts that by increasing $N$  one can make $p_k$ arbitrary close to $\widehat p_k$ (uniformly).
We need to estimate the length of a zig-zag. Summing up lengths of segments of the orbit of $p$ and adding and subtracting $\sum_{k=0}^{N-1} \sum_{i=1}^r \length (e^{\frac{t}{N}X_i}p_k)$ we obtain
\begin{align}
\length((e^{\frac{t}{N}X_1}\dots e^{\frac{t}{N}X_r})^Np)=\sum_{k=0}^{N-1}\sum_{i=1}^r
\length (e^{\frac{t}{N}X_i}p_k)+ \label{eq:1}\\
+ \sum_{k=0}^{N-1} \Big(\length (e^{\frac{t}{N}X_1}\dots e^{\frac{t}{N}X_r}\widehat p_k) - \sum_{i=1}^r \length (e^{\frac{t}{N}X_i}p_k)\Big)\label{eq:2}
\end{align}

Figure \ref{fig:zig-zag} illustrates the case of two vector fields $X$ and $Y$. By Lemma \ref{abraham} again and continuity of flows, for any $\epsilon>0$ we can chose $N$ big enough such that for any $i$ and all $k$  one has
$\big\lvert e^{\frac{t}{N}X_{i+1}}\dots e^{\frac{t}{N}X_r} \widehat p_k - p_k \big\rvert<\epsilon$.
Applying Lemma \ref{stretch} for $X=X_i$ one obtains for each $i$:

$$\big\lvert\length (e^{\frac{t}{N}X_i}\dots e^{\frac{t}{N}X_r}\widehat p_k) -  \length (e^{\frac{t}{N}X_i}p_k)\big\rvert<C\epsilon\frac{t}{N}$$

Summing up over $i$ and applying the triangle inequality for the absolute value and the additivity of length, we obtain

$$
\big\lvert \length (e^{\frac{t}{N}X_1}\dots e^{\frac{t}{N}X_r}\widehat p_k) - \sum_{i=1}^r \length (e^{\frac{t}{N}X_i}p_k)\big\rvert<
Cr\epsilon\frac{t}{N}
$$

And now summing up over $k$ and applying the triangle inequality again

$$
\big\lvert \sum_{k=0}^{N-1}\Big(\length (e^{\frac{t}{N}X_1}\dots e^{\frac{t}{N}X_r}\widehat p_k) - \sum_{i=1}^r \length (e^{\frac{t}{N}X_i}p_k)\Big)\big\rvert<
Cr\epsilon t
$$

Thus the sum (\ref{eq:2}) in the expression for the length of a zig-zag tends to zero when $\epsilon\to 0$ (while $N\to \infty$). On the other hand, by the definition of the length of the path we can differentiate it:

$$\frac{d}{ds}\length(\gamma(s))\Big|_{s=0}=\norm{\dot\gamma(0)}$$

For an integral curve emanating from a point $p$ this means that for any $\frac{\epsilon}{r}>0$, any vector field $X_i$ and $N$ large enough
$$\Big\lvert\frac{N}{t}\length(e^{\frac{t}{N}X_i}p_k)-\norm{X_i(p_k)}\Big\rvert<\frac{\epsilon}{r}$$

Hence, summing up, we obtain
\begin{gather*}
\Big\lvert\sum_{k=0}^{N-1}\sum_{i=1}^r
\length (e^{\frac{t}{N}X_i}p_k) - \sum_{k=0}^{N-1}\sum_{i=1}^r
\frac{t}{N}\norm{X_i(p_k)} \Big\rvert< \\ < \sum_{k=0}^{N-1}\frac{t}{N}\epsilon= \epsilon t
\end{gather*}

However $\sum_{k=0}^{N-1}\frac{t}{N}\sum_{i=1}^r
\norm{X_i(p_k)}$ is just a Riemann sum for the integral $\int_0^t \sum_{i=1}^r \norm{X_i(\gamma(s))}ds$ of the function $s\mapsto \sum_{i=1}^r (\norm{X_i(\gamma(s))}$.
The function $\sum_{i=1}^r\norm{X_i(\gamma(s))}$ is continuous as a function of $s$ and, consequently, Riemann integrable. Thus $$\lim_{N\to\infty}\sum_{k=0}^{N-1}
\frac{t}{N}\sum_{i=1}^r\norm{X_i(p_k)} = \int_0^t \sum_{i=1}^r(\norm{X_i(\gamma(s))}ds$$. Hence the sum (\ref{eq:1}) converges: 

$$
\lim_{N\to \infty} \sum_{k=0}^{N-1}\sum_{i=1}^r
\length (e^{\frac{t}{N}X_i}p_k)=\int_0^t \sum_{i=1}^r\norm{X_i(\gamma(s))}ds
$$

\end{proof}

\begin{center}
\begin{minipage}{0.7\textwidth}
  \centering
  \includegraphics[width=1\linewidth]{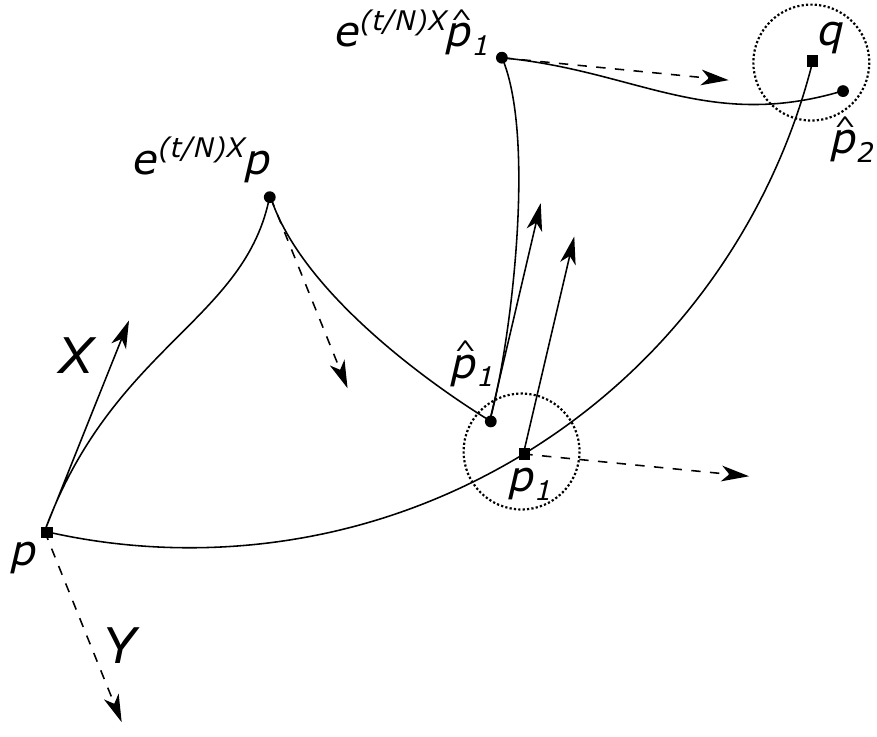}
  \captionof{figure}{Zig-zag approximation}
  \label{fig:zig-zag}
\end{minipage}%
\end{center}

\end{subsection}

Note that if all $X_i$ lie in the distribution of cones $D$, the corresponding zig-zags are $D$-admissible paths by definition. We are going to use them to approximate all $\langle D \rangle$-admissible paths. We will denote a zig-zag path $$\big(e^{\frac{t}{N}X_1}\cdots e^{\frac{t}{N}X_{m}}\big)^Np$$ by $Z_\mathbf{X}^N(p,t)$, where $\mathbf{X}=(X_1,\dots, X_m)$.

Next we need to show that the Finsler (i.e. measured with respect to the norm $\norm{\cdot}_{g,D}$) length of a $\langle D\rangle$-admissible path is close to the Riemannian length of some zig-zag path.

\begin{lemma}\label{norm_realization}
Let $\gamma$ be a  $\langle D\rangle$-admissible path in $U$. For any $\epsilon>0$ there exists a finite set of vector fields $X_1, X_2,\dots X_m$ from $D$, such that $\dot\gamma=X_1 + X_2+\dots + X_m$ and at every point along the path
\begin{equation}\label{realiz}
\big\lvert|X_1|+|X_2|+\dots+ |X_m| - \norm{\dot\gamma}_{g,D}\big\rvert<\epsilon
\end{equation}

\end{lemma}

\begin{remark}
In other words, we can always choose global vector fields from $D$ that approximately realize the norm $\norm{\cdot}_{g,D}$ at each point. Since we do not require the cone $D_p$ to be closed (only the fiberwise closure of $D$ to be closed, as in Remark \ref{fiberwise_clos}), the norm need not be realizable exactly by a set of vectors even infinitesimally.
\end{remark}

\begin{proof}

As in the proof of Proposition \ref{cont_norm},  in every tangent space $T_pM$ along $\gamma$ we can decompose $\dot\gamma$ as a finite sum $\sum_i X_{ip}$ of linearly independent vectors from $D_p$ that satisfy the inequality \eqref{realiz} with $\epsilon/4$ on the right-hand side. Each of those vectors by Definition \ref{conedef} is a value of some global vector field $\overline X_{ip}$ from the cone $D$. Thus for each point $p$ by continuity we can choose a small neighborhood $U_p$ such that inequality (\ref{realiz}) is satisfied with $3\epsilon/4$ on the right-hand side. If $\overline X_{ip}$ are linearly independent at $p$ they are also linearly independent in some neighborhood of $p$. The vector fields $X'_{ip}:=\frac{\langle \dot{\Gamma}, \overline X_{ip} \rangle}{|\dot\Gamma|}\overline X_{ip}$ (where $\Gamma$ is an arbitrary extension of $\dot\gamma$ to a vector field in a neighborhood of $\gamma$) decompose $\dot\gamma$, that is $\dot\gamma=\sum_i X'_{ip}$ on each $U_p$. By shrinking further the neighborhood $U_p$, if necessary, we can ensure that the inequality (\ref{realiz}) holds for $X'_{ip}$ in $U_p$.

Now choose a finite number of $U_p$ that cover $\gamma$. Taking a partition of unity $\{\phi_p\}$ subordinate to $U_p$, define $X_i:=\sum_p \phi_p  X'_{ip}$. At each point one has $\dot{\gamma}=\sum_i X'_{ip}$ and hence $\dot{\gamma}=\sum_i X_i$. By the homogeneity of the norms $\norm{\cdot}_{g,D}$ and $|\cdot|$ the inequality (\ref{realiz}) is also satisfied.

\end{proof}

Thus, combining Lemmas \ref{norm_realization} and \ref{mainlemma}, we obtain the statement

\begin{proposition}
\label{zig-zag}
For any $\langle D\rangle$-admissible path $\gamma$ between points $p,q\in M$ and any $\epsilon>0$ there exists a $D$-admissible path $\widehat\gamma=Z_\mathbf{X}^N(p,t)$ between the points $p$ and $q'$, such that $\big\lvert q-q' \big\rvert<\epsilon$ and $\big\lvert\length( \widehat\gamma) - \length_{sF}(\gamma)\big\rvert<\epsilon$
\end{proposition}

Note that the approximating path does not connect the same points $p$ and $q$. However, if we can connect $q'$ with $q$ by a $D$-admissible path $\delta$ of length less than $\epsilon$, then the Theorem \ref{mainthm} follows. Indeed, the concatenation of a zig-zag $\widehat\gamma$ and $\delta$ will be $D$-admissible, connecting $p$ and $q$ and deviating from $\length_{sF}(\gamma)$ not more than by $2\epsilon$. The connecting of $q'$ with $q$ is the subject of the next subsection.

\subsection{Accessible set}
The proposition below is essentially the Chow-Rashevski theorem (or Sussmann orbit theorem). It is a basic result of subriemannian geometry and geometric control theory \cite{gromov},\cite{Agrachev},\cite{bellaiche}. However, it is not always stated in the form and generality that is required for our situation. Hence, for the sake of continuity of the exposition, we sketch the proof here.

Let $\C$ be a finite set of vector fields on an open neighborhood of $0$ in $\mathbb{R}^n$. Let us denote the set of all possible Lie brackets of vector fields from $\C$ by $\Lie(\C)$. Suppose $\Lie(\C)$ span $T\mathbb{R}^n$. Let us define a finite set of flows $\Phi_i^t$ as follows:
\begin{enumerate}
    \item Each of the flows $e^{tX}$ of vector fields $X\in \C$ are among $\Phi_i^t$.
    \item For every pair $X,Y$ of fields $X,Y\in \C$ define a flow $$\Phi^{t}_{[X,Y]}:=e^{-\sqrt tY}e^{-\sqrt tX}e^{\sqrt tY} e^{\sqrt tX}$$
    \item For higher bracket expressions of vector fields from $\C$ by induction, assuming we already defined it for bracket expressions $A,B$ of smaller length:
    $$\Phi^{t}_{[A,B]}:=\Phi^{-\sqrt[l]t}_A \Phi^{-\sqrt[l]t}_B\Phi^{\sqrt[l]t}_A \Phi^{\sqrt[l]t}_B$$
    where $l$ is a length of the bracket expression $[A,B]$.

\end{enumerate}

Let us enumerate $\Lie(\C)$ and hence all such flows: $\Phi_i^t$, $i=1,\dots m$. Define the end point map as follows:
$$E: (t_1,t_2,\dots t_m) \mapsto \Phi_1^{t_1}\Phi_2^{t_2}\dots \Phi_m^{t_m}$$

\begin{lemma}
\label{chow}
The end point map $E$ is open at $(0,\dots 0)$.
\end{lemma}

\begin{proof}
The end-point map $E$ is a smooth map from an open neighborhood of $0$ in $\mathbb{R}^m$ to $\mathbb{R}^n$. A higher commutator of flows of vector fields is tangent to the corresponding Lie bracket of the fields to the order $l$ (hence power of $l$ in the definition of the flows) \cite[Theorem 1]{michor}. It follows that that $E$ is surjective at $(0,\dots 0)$. By the constant rank theorem the image of a small enough open neighborhood of $(0,\dots 0)$ is open.
\end{proof}
\begin{remark}
Note that each $\Phi_i^t(p)$ is connected with $p$ by a $\C$-admissible path by construction. For a nice explicit description in the three-dimensional case see \cite[Theorem 5.4.6]{BBI}.
\end{remark} 

\begin{proposition}
Let $\C$ be an arbitrary distribution of cones. For every $\delta>0$ there is an open neighborhood $U$ of $0$ in $\mathbb{R}^n$ such that  any point of $U$ can be connected to $0$ by a $\C$-admissible path of length less than $\delta$.
\end{proposition}
\begin{proof}
Chose a finite set of vector fields from $\C$ satisfying the Chow condition. By the Lemma \ref{chow} for any small $\delta>0$, the image of a $\delta$-ball around $0\in \mathbb{R}^m$ is an open neighborhood $U$ of $0\in \mathbb{R}^n$.

The length of the a path $\Phi_1^{t_1}\Phi_2^{t_2}\dots \Phi_m^{t_m}(0)$ corresponding to $(t_1,t_2,\dots t_m)$ can be estimated from above. Let $K$ be the supremum of norms of vector fields from $Lie(\C)$ on some compact neighborhood of $0\in\mathbb{R}^n$. Then 
$$\length (\Phi_1^{t_1}\Phi_2^{t_2}\dots \Phi_m^{t_m}(0))\le K(|t_1|+|t_2|+\dots+|t_m|)$$

Choosing $\delta$ small enough, we obtain an open neighborhood of $0\in\mathbb{R}^n$, each point of which is connected to $0$ by a $\C$-admissible path $\Phi_1^{t_1}\Phi_2^{t_2}\dots \Phi_m^{t_m}(0)$ of length less than $\delta$.
\end{proof}

\section{Sub-twistor metric}
\subsection{Homogeneous cones}
A family of submanifolds in a manifold can give rise to a natural sub-conical metric, which can be thought of as a higher-dimensional generalization of the Manhattan distance, where streets are higher dimensional submanifolds.

\begin{example} Let $S^2_+$ be the subset of the unit two-dimensional sphere in $\mathbb{R}^3$, consisting  of all points with $z$-coordinate between $-\sqrt 2$ and $\sqrt 2$. At any point $x\in S^2_+$ the set of vectors tangent to great circles that pass through $x$ and {\it lie entirely in $S^2_+$} sweep out an open cone that we define to be the cone of admissible directions $D_{sph}$. Consider the corresponding sub-conic metric.

The shortest path connecting two points in general may have ``breaks'', since  a great circle connecting two points may not lie entirely in $S^2_+$ (Figure 2).
\end{example}

\begin{center}
\begin{minipage}{.5\textwidth}
  \centering
  \includegraphics[width=1\linewidth]{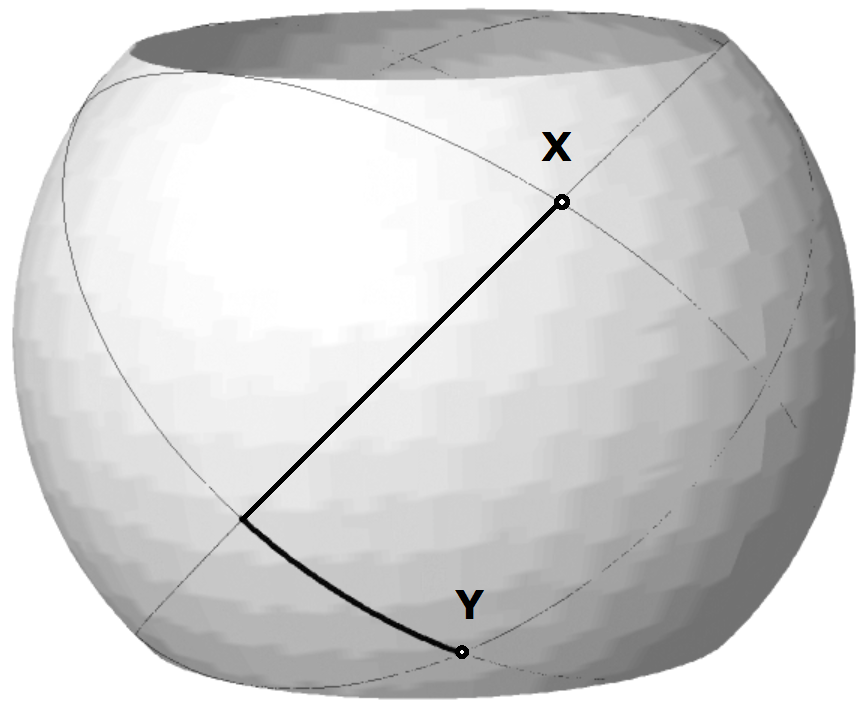}
  \captionof{figure}{Positive sphere metric}
  \label{fig:sphere}
\end{minipage}
\end{center}

\bigskip

The metric space $S^2_+$ is in a sense a toy model of the situation we will encounter later.
\begin{defn}
\label{submanifolds}\label{tranhom}
Let $M=G/H$ be a homogeneous manifold and $\mathcal{L}$ a family of submanifolds, such that 
\begin{enumerate}
\item For any point $x\in M$, there is a submanifold from $\mathcal{L}$ passing through $x$. Denote the set of such manifolds by $\mathcal{L}_x$

\item For any two points $x,y\in M$ there is a finite sequence of submanifolds $S_0\dots S_n$ from $\mathcal{L}$, such that $S_i\cap S_{i+1}\neq \varnothing$ and $x\in S_0$, $y\in S_n$.

\item The group maps submanifolds from $\mathcal{L}$ to submanifolds from $\mathcal{L}$.
\end{enumerate}
We will call such a family of submanifolds $\mathcal L$ {\it transitive homogeneous}.
\end{defn}

\begin{proposition}
\label{rich}
The union of all tangent spaces  of submanifolds from a transitive homogeneous family $\mathcal L$ (viewed as subsets of $TM$) is a distribution of cones. We will denote it by $\mathcal{C(L)}$.
\end{proposition}

\begin{proof}
The first condition of Definition \ref{tranhom} rules out pathological cases with isolated points. The second condition ensures finitness of the resulting metric ruling out  foliation-like behavior. 

To show that $\mathcal{C(L)}$ is a distribution of cones we need to demonstrate that each vector tangent to a submanifold from $\mathcal{L}$ at some $p\in M$ extends to a vector field $X$ in a neighborhood of $p$ such that $X_p$ is tangent to some submanifold from the family $\mathcal{L}$. 

Let $v$ be a vector from $T_{p_0}S\subset T_{p_0}M$, where $S\in \mathcal{L}$. 
Since $M=G/H$ is a homogeneous manifold for the group $G$, then $\operatorname{pr}: G\to G/H$ is a locally trivial bundle. Choose a smooth section $\nu$ of $\operatorname{pr}$ in a neighborhood $U$ of $p_0\in M$ such that $\nu(p)=1\in G$. The vector field $p\mapsto \nu(p)_*v$ is defined in an neighborhood of $p_0$, smooth and is tangent at a point $\nu(p)(p)$ to a submanifold $\nu(p)(S)$ which belongs to $\mathcal{L}$ by the third condition of the proposition.

It is left to prove that the fiber-wise closure $\overline{\mathcal{C(L)}}$ of $\mathcal{C(L)}$ is closed. Suppose the sequence of tangent vectors $v_i\in \overline{\mathcal{C(L)}}$ converges to some vector $v\in T_pM$ and $v\notin \overline{\mathcal{C(L)}}$. By homogeneity and Property 3 of Definition \ref{tranhom}, for each $i$  pick an element of the group $g_i\in G$ such that $g_iv_i\in T_pM$  and $g_i\to Id$ as $i\to\infty$. Each $g_iv_i\in T_pM$ lies in a closed subset (the closure of the cone in the fiber) and converges to $v$, which yields a contradiction.

\end{proof}

\subsection{The period space}
This subsection is devoted to the main subject of the paper --- the period space $\per$. It is introduced here as a purely linear-algebraic object. Its meaning as a period space for a hyperk\"ahler manifold is explained in the appendix.

Let $V$ be a vector space over $\mathbb{R}$ of dimension $b$ and $q$  a quadratic form on $V$ of signature $(3,b-3)$. We will write $V>0$ to indicate that $V$ is a positive-definite subspace.
\begin{defn}\label{perdefn}
{\it The period space} is the set of all oriented positive definite two dimensional subspaces in $V$ :
$$\per:=\left\{E:E\subset\text{V oriented subspace of } \dim E=2\text{ and } q(E)>0\right\}$$
\end{defn}

\begin{remark}
Period space is easily seen to be a homogeneous space for the group of orientation preserving orthogonal transformations for the metric of signature $(3,b-3)$:
$$\per\simeq SO(3,b-3)/SO(2)\times SO(1,b-3)$$
\end{remark}

Period space is connected and simply connected \cite[Claim 2.9]{Ver}.

\begin{remark}
$\per$ is an open subset of the space $Gr^{+}(2,V)$ of all oriented two-dimensional subspaces of $V$. In particular we can identify the tangent space of $\per$ at a point with the tangent space to $Gr^{+}(2,V)$. That is, if a point of $\per$ represented by a plane $W$ then $T\per$ can be canonically identified with $\mathrm{Hom}(W,W^\perp)$.
\end{remark}
Now we are ready to define an important family of submanifolds in $\per$:
\begin{defn}
Let $T$ be a positive $3$-dimensional subspace of $V$. The {\it twistor sphere}, associated to $T$ is the set of all oriented planes in $T$.
$$\Tw_T:=\big\{L: L\in Gr^+(2,V) \text{ and } L\subset T\big\}$$
\end{defn}
\begin{remark}
A twistor sphere is homeomorphic to a 2-dimensional sphere.
\end{remark}

\begin{defn}Tangent spaces of all twistor spheres passing through a point $p\in\per$ sweep out a cone which we will call a {\it subtwistor cone}  $D^{tw}$. In other words, for $p\in \per$
$$D^{tw}_p:=\bigcup_S T_pS$$
where union is taken over all twistor spheres $S$
 passing through $p$. Here we view every $T_p S$ as a subspace of $T_p \per$.

\end{defn}
In the next subsection we will describe the shape of this cone at a tangent space at a point and show that it is a distribution of cones.
\subsection{Sub-twistor cone}

\begin{proposition}
\label{twistor_tangent}
At a point $W\in \per$ the subtwistor cone is the subset of $\Hom(W,W^\perp)$ defined by the following conditions:
\begin{enumerate}
    \item Rank of $v\in \Hom(W,W^\perp)$ does not exceed one,
    \item Image of $v$ is positive with respect to the restriction of the form $q$ on $W^\perp$, or $v=0$
\end{enumerate}
\end{proposition}

\begin{proof}
Actually the second condition implies the first since the signature of the form $q$ restricted to $W^\perp$ is $(1,b-3)$ and there are no positive subspaces of dimension greater than $1$ in it.

Given a positive plane $W$ we can choose a positive vector  $v$ orthogonal to $W$ and take the linear span $T=\langle W,v\rangle$. Any positive $3$-subspace containing $W$ is  uniquely defined by such $v$. 

As usual, we can identify planes in an open neighborhood of $W$  with the graphs of linear maps $\Hom(W,W^\perp)$. If a tangent vector $u\in T_W\per$ is tangent to the twistor sphere $\Tw_T$ then its graph is a subset of $T$. Consequently, the image of $u\in \Hom(W,W^\perp)$ is $\mathbb{R}v$. Thus for any $v$ in the positive cone of $W^\perp$ there is a $2$-dimensional set $\Hom(W,\mathbb{R}v)$ of tangent vectors to $\per$ tangent to the corresponding twistor sphere $\Tw_T$. The union of all such sets is the subtwistor cone:

$$D_W^{tw}=\bigcup_{\substack{L>0 \\ L\subset W^\perp  }} \Hom(W,L)$$

\end{proof}

\begin{remark}

Now we can describe the sub-twistor cone in coordinates by a system of  equations and inequalities. Let us choose a basis $e_1,e_2,\dots e_b$ in $V$ such that the form $q$ is diagonal $(+++,--\dots)$ and $W=\langle e_1,e_2\rangle$. A $(2,b-2)$ matrix $A$ is zero or has a positive definite image if and only if the following two condition hold:
\begin{enumerate}
    \item Its rank is not greater than  $1$. In other words, $A(e_1)$ is proportional to $A(e_2)$.
    \item $A(e_1)=A(e_2)=0$ or $q(A(e_1))+q(A(e_2)) > 0$. Indeed, $q(A(e_1))$ and $q(A(e_1))$ are both non-positive or non-negative and we want at least one of them to be strictly positive unless $A$ is zero.
\end{enumerate}

The first condition is equivalent to the vanishing of all $(2,2)$ minors of $A$. That is, it is a system of $b-3$ quadratic equations in $2b-4$ variables. The second condition is equivalent to a quadratic inequality in $2b-4$ variables. 

\end{remark}

\begin{proposition}
Subtwistor cone $D^{tw}$ is a distribution of cones on $\per$. 

\begin{proof}

We saw that there are many twistor spheres passing through each point of $\per$. Any two points can be connected by a sequence of twistor spheres. This is proved in \cite[Proposition 3.7]{bourb} and \cite[Proposition 5.8]{Ver}.

The period space is a homogeneous space for the group $G=SO(3,b-3)$. $G$ preserves the form $q$ and thus acts on the set of $3$-dimensonal positive subspaces of $V$. Hence elements of $SO(3,b-3)$ on $\per$ map twistor spheres to twistor spheres.

 Thus set of all twistor spheres is a transitive homogeneous family of submanifolds (Definition \ref{tranhom}). Hence by Proposition \ref{rich} it is a distribution of cones.
\end{proof}

\end{proposition}

Thus the sub-twistor metric satisfies the conditions of Theorem \ref{mainthm}. The next proposition describes the linear hull of the sub-twistor cone.
\begin{proposition}
The linear hull of the sub-twistor cone is the whole tangent bundle, that is $\langle D^{tw} \rangle = T\per$. In particular, the Chow condition is trivially satisfied.
\end{proposition}

\begin{proof}
An arbitrary  $2\times (b-2)$ -matrix can be written as a sum of two matrices of rank $1$. Thus to prove the proposition we need to be able to express a rank $1$ linear map as a sum of matrices with positive images. 

Let $A$ be an arbitrary matrix of rank $1$. Pick a vector $v\in W$, $v\notin Ker(A)$.
Let $B$ be another matrix with the same kernel $Ker(B)=Ker(A)$ and for which $q(B(v))\gg 0$. Then $q(A(v)+B(v))>0$ hence $A+B$ has  positive image. But $-B$ is also positive on $v$ by symmetry of $q$ and we have
$$A=(A+B)-B$$
So any rank $1$ matrix can be written as a sum of two matrices with positive image, that is, matrices from $D^{tw}$. Thus any tangent vector from $T\per$ can be written as a sum of at most four vectors from $D^{tw}$

\end{proof}

\subsection{Sub-twistor norm}
Applying Theorem \ref{mainthm} we obtain the following

\begin{theorem}
\label{periodfinsler}
Sub-twistor metric is Finsler. 
\end{theorem}
The norm is defined by Theorem \ref{mainthm} as follows: intersect the sub-twistor cone $D^{tw}$ with the unit riemannian ball $B$ and take the closure of  the convex hull of $D^{tw}\cap B$. The resulting convex body will define the norm by the Minkowski construction. 

\begin{remark}
Note, that this convex hull of $D^{tw}\cap B$ is not strictly convex in the tangent space at a point. That is, there is a hyperplane  which intersects the closure  of $D^{tw}\cap B$ (and is disjoint with interior) in at least two points (hence the convex hull will contain the interval connecting these two points). The affine hyperplanes intersecting the unit ball are in one to one correspondence with boundaries of geodesic balls on $\partial B$ (given by intersection of a plane with  $\partial B$), so the claim can be deduced from the following proposition (cf. \cite{ros}):

\begin{proposition}
Let $C$ be a compact proper subset of an $n$-sphere $S^n$. Then there exists a geodesic ball $K$ in the complement of $C$, such that $\partial K\cap \partial C$ contains at least two points, and $K$ and $C$ have disjoint interiors.
\end{proposition}

\begin{proof}
The space of geodesic balls in the closure of the complement of  $C$ is compact. Choose a ball $K$ of maximal radius $r$. Suppose it intersects $\partial C$ in a point $p$. Divide $K$ by a geodesic hyperplane into two congruent halves, one of which contains $p$. Let us denote another one by $H$. We claim that $\partial K$ also intersects $\partial C$ at some point of $H\cap \partial K$. Otherwise, by compactness of $H\cap \partial C$ and $C$, there exists a positive number $\varepsilon$ such that $\varepsilon$-neighborhood of $H$ does not intersect $C$.
The statement now follows from the following geometric property of spherical (as well as euclidean) balls: the union of $K$ and the $\varepsilon$-neighborhood of $H$ contains a ball of radius bigger than $r$.

\end{proof}
Since the sub-twistor norm is not strictly convex then the minimizing paths are not locally unique \cite[Section 5.1.2]{BBI}.
\end{remark}

\begin{remark}
Finsler nature of sub-twistor metrics implies that the induced topology coincides with the topology induced by the ambient riemannian metric and hence with the underlying manifold topology. This is the crucial property in the application of sub-twistor metrics to the proof of the global Torelli theorem. In \cite{Ver}  the solution of Gleason and Palais to a generalization of the Hilbert fifth problem  was used to obtain this fact.
\end{remark}

\begin{remark}
Note that any vector $v$ from the sub-twistor cone at $W$ is tangent to a single twistor sphere. Namely, the twistor sphere associated to a positive 3-space spanned by the image of $v$ in $W^\perp$ and $W$.
\end{remark}

\begin{proposition}
Let $W_0\in \per$ and $\Tw_{T_0}$ be a twistor sphere passing through $W$. There exists a neighborhood $U$ of $W_0$ and a foliation $\mathcal F$ of $U$ by intersections of twistor spheres with $U$ such that $\Tw_{T_0}$ is the leaf through $W_0$. Also, any tangent vector $w\in T_{W_0}$ can be extended locally to a vector field tangent to $\mathcal F$.
\end{proposition}
\begin{proof}
A tangent vector $w$ to a twistor sphere $\Tw_{T_0}$ passing through $W_0$ defines a positive vector $v\in W_0^\perp$,unique up to multiplication by a constant (see Proposition \ref{twistor_tangent}). Now consider the collection $\mathcal{F}$ of all positive $3$-planes $T$ containing $v$. We claim that $\mathcal{F}$ is a foliation of an open neighborhood of $W_0$. 

\begin{enumerate}[Step 1:]
\item Since positivity is an open condition, there is an open neighborhood of $W_0$ consisting of positive $2$-planes $W$ such that the space generated by $W$ and $v $ is positive. Hence there is a twistor sphere  from $\mathcal{F}$ passing through each point of this neighborhood.
\item The set of all positive $2$-planes that do not contain $v$ is open, intersect it with the neighborhood from the previous step. Now suppose two positive $3$-subspaces containing $v$ have intersection of dimension $2$ which does not contain $ v $. Then their intersection necessarily has dimension $3$, so they coincide.
\item Now we can extend the given tangent vector $w$ to a vector field in some neighborhood which is tangent to $\mathcal F$. Indeed,  trivializing the tangent bundle in the neighborhood of $W_0$ we can identify $T_p\per$ with $\Hom(W_0,W_0^\perp)$. Now pick any extending vector field (Proposition \ref{rich}) and compose it with a projection on $v$. This vector field will be tangent to some twistor sphere from $\mathcal F$ at each point.
\end{enumerate}

\end{proof}

\begin{remark}
\label{piecewise_spheric}
By the theorem on the uniqueness of solutions of ordinary differential equations the trajectory of any point of $U$ under the flow of a vector field of the proposition above lie on a twistor sphere. Thus we can ensure that zig-zags from the proof of the main Theorem \ref{mainthm} are piecewise smooth curves, each smooth segment of which lying on one twistor sphere. Thus the cone defining the sub-twistor metric can be assumed to consist of vector fields locally tangent to foliations by twistor spheres.

\end{remark}

\appendix
\section{Global Torelli theorem}

In this appendix we give an overview of Verbitsky's global Torelli theorem for hyperk\"ahler manifolds. This gives the motivation and context for the introduction of the space $\per$ and the sub-twistor metric on it. 

Broadly speaking, Torelli-type theorems aim to describe moduli spaces of geometric structures in terms of Hodge-theoretic data. 

\subsection{Hyperk\"ahler manifolds}

\begin{defn}
A hyperk\"ahler manifold is a Riemannian manifold $M$ equipped with a triple of complex structures $I,J$ and $K$ such that:

\begin{enumerate}
    \item $I^2=J^2=K^2=IJK=-1$
    \item $M$ is K\"ahler with respect to each complex structure $I,J$ and $K$. We will denote the corresponding symplectic forms by $\omega_I, \omega_J, \omega_K$
\end{enumerate}

In this paper we {\it always assume} that a hyperk\"ahler manifold is simply connected and compact.
\end{defn}

\begin{remark}
The form $\omega_J+\sqrt{-1}\omega_K$ is a complex valued symplectic form which is holomorphic with respect to the complex structure $I$, that is, it is of type $(2,0)$.
\end{remark}

\begin{defn}
A compact simply connected K\"ahler manifold that possesses a  holomorphic symplectic form with $H^{2,0}(M,\mathbb{C})$  one-dimensional is called an irreducible holomorphic symplectic manifold (IHS for short).
\end{defn}

As a consequence of the celebrated Calabi-Yau theorem one obtains the following fact: every compact K\"ahler manifold admitting a non-degenerate holomorphic $2$-form is a hyperk\"ahler manifold.

\begin{proposition}
The bundle automorphism
$$L=a I+ b J+ c K$$
where $a^2+b^2+c^2=1$ is an almost complex structure, that is, $L^2=-1$. Moreover, $L$ is integrable and K\"ahler with respect to the Riemannian metric on $M$.
\end{proposition}
In other words, given a hyperk\"ahler structure we obtain a family of K\"ahler  structures parametrized by the  two dimensional sphere.

The second cohomology  $H^2(M,\mathbb{Z})$ of a hyperk\"ahler manifold admits a remarkable quadratic form which generalizes the intersection form in real dimension $4$.
\begin{proposition}\label{bbf} Let $M$ be an IHS manifold of real dimension $4n$. Then there exists a primitive quadratic form $q: H^2(M,\mathbb{Z})\to \mathbb{Z}$  with the following property: there exists a rational constant $c>0$, such that for any class $\eta \in H^2(M,\mathbb{Z})$

$$\int_M \eta^{2n}=c\cdot q(\eta,\eta)^n$$
\end{proposition}
This form is called {\it Beauville-Bogomolov-Fujiki form} and it is the unique form satisfying Proposition \ref{bbf} up to sign. The sign is chosen in such a way that the following formula holds (for some $\lambda>0$):

\begin{gather*}
    \lambda q(\eta,\eta')=2\int_M \eta\wedge\eta'\wedge \Omega^{n-1}\wedge \bar\Omega^{n-1}-\\
    \frac{n}{n-1}\frac{(\int_M \eta\wedge \Omega^{n-1}\wedge\bar\Omega^{n})(\int_M \eta'\wedge \Omega^{n}\wedge\bar\Omega^{n-1})}{\int_M \Omega^n\wedge\bar\Omega^n}
\end{gather*}

\begin{proposition}
The signature of the Beauville-Bogomolov-Fujiki form is $(3,b-3)$ where $b=\dim H^2(M)$. A $3$-space  spanned by the holomorphic symplectic form $\Omega$, its complex conjugate $\bar\Omega$ and a K\"ahler form $\omega$ is positive.
\end{proposition}

\begin{remark}\label{quadric}
$q(\Omega,\bar\Omega)>0$ and $q(\Omega,\Omega)=0$.
\end{remark}

\subsection{Teichm\"uller and birational Teichm\"uller spaces}

\begin{defn}
We say that a a complex structure $(M,I)$ is of {\it hyperk\"ahler type} if it is a complex structure underlying some IHS manifold, that is, simply connected compact with a holomorphic (with respect to the complex structure $I$) symplectic form that spans $H^{2,0}(M)$.
\end{defn}

Consider the space of all complex structures on $M$. That is, the set of sections $J$ of the bundle of tensors $\End(TM)$ satisfying the equation $J^2=-1$ and the integrability condition. 

Now let $\comp$ be a subspace of complex structures of hyperk\"ahler type on $M$, equipped with the topology of uniform convergence of all derivatives, and $\texttt{Diff}_0$ the group of diffeomorphisms of $M$ isotopic to the identity. Diffeomorphisms from $\texttt{Diff}_0$ act by pullbacks on complex structures of hyperk\"ahler type. 
\begin{defn}
The quotient space $\comp/\diff$ is called the {\it Teichm\"uller space} of $M$ and denoted by $\teich$. 
\end{defn}

The results of Kuranishi and the Bogomolov-Todorov-Tian theorem imply that $\teich$ is a complex manifold (see \cite[Remark 1.7]{Ver}). However, unlike in the situation of Riemann surfaces, it is generally not Hausdorff. To formulate global Torelli theorem one has to work with a so called Hausdorff reduction of $\teich$.

\begin{defn}
Let $X$ be a topological space. A {\it Hausdorff reduction} of $X$ is a continuous surjection $r: X\to \bar{X}$, such that any continuous map to a Hausdorff space $f:X\to Y$ factorizes through $r$.
\end{defn}
A Hausdorff reduction need not exist in general. 
\begin{defn}
Let $p$ and $q$ be points in a non-Hausdorff topological space $X$. They are called inseparable (denoted by $p\sim q$) if for every open neighborhood $U$ of $p$ and open neighborhood $V$ of $q$ we have $U\cap V\neq \empty \varnothing$.
\end{defn}

It may be shown \cite[Section 4]{Ver} that $\teich/\sim$ is a Hausdorff reduction for $\teich$. The crucial ingredient in the proof of this fact is the following theorem of Huybrechts \cite{Huy}:

\begin{theorem}\label{huybirat}
Let $I,J\in \teich$ and they are inseparable. Then $(M,I)$ and $(M,J)$ are bimeromorphic.

\end{theorem}

\begin{defn}
The birational Teichm\"uller space $\teich_b$ is the space $\teich/\sim$. It is a Hausdorff reduction of $\teich$.
\end{defn}

It is a Hausdorff topological space. Moreover, the local Torelli theorem of Bogomolov \cite{bog} implies that it is actually a manifold.

\subsection{Period map and period domain}

Let $(M,I)$ be an IHS manifold and denote by $H^{p,q}(M,I)$ the Hodge decomposition of the cohomology of $M$.

\begin{proposition}
\label{diff_cohom}
If $J$ is a complex structure obtained as a pullback of $I$ by a diffeomorphism $f\in \texttt{Diff}_0$ then $H^{p,q}(M,I)=H^{p,q}(M,J)$.
\end{proposition}

In other words, the Hodge decomposition is well-defined on the Teichm\"uller space \cite[Claim 1.9]{teich}. The period map associates to an IHS manifold its $(2,0)$-Hodge component $H^{2,0}(M)$:

\begin{defn}
The period map $\Pi: \teich \to \mathbb{P}H^2(M,\mathbb{C})$ is defined by:

$$\Pi: I\mapsto H^{2,0}(M,I)\subset H^2(M,\mathbb{C})$$
\end{defn}
This map does not depend on the choice of a representative of a complex structure in $\teich$ by  Proposition \ref{diff_cohom}. Since all complex structures $I$ are of hyperk\"ahler type then corresponding $H^{2,0}(M,I)$ are one dimensional complex subspaces and hence are indeed elements of $\mathbb{P}H^2(M,\mathbb{C})$.

It is not hard to see that $\Pi$ factorizes through the Hausdorff reduction $\teich\to\teich_b$. By Theorem \ref{huybirat} different preimages in $\teich$ of $I\in \teich_b$ are bimeromorphic. Clearly, such manifolds $M$ have the same spaces of holomorphic forms. From now on we will be referring to the map $$\Pi: \teich_b\to \mathbb{P}H^2(M,\mathbb{C})$$ as the {\it period map}.

\begin{defn}
Note that by Remark \ref{quadric} image of $\Pi$ is contained in the set defined by 
    \[   \left\{
\begin{array}{ll}
      q(\Omega,\Omega)=0 \\
      q(\Omega,\bar\Omega)>0
\end{array} 
\right . \] 
We will call it {\it the period domain}.
\end{defn}

Identification of the period domain  with the period space $\per$ in the sense of Definition \ref{perdefn} is given by the following proposition \cite[Proposition 4.8]{teich}:

\begin{proposition}
The space $\{\Omega\in H^2(M,\mathbb{C}): q(\Omega,\Omega)=0, q(\Omega,\bar\Omega)>0\}$ and the Grassmannian of positive oriented $2$-dimensional subspaces in $H^2(M,\mathbb{R})$ are diffeomorophic.
\end{proposition}

Now we are ready to state the global Torelli theorem:

\begin{theorem}
The map $\Pi: \teich_b\to \per$ is a diffeomorphism on connected components of $\teich_b$.
\end{theorem}

\begin{proof}
This is a direct consequence of the fact that the period map is a covering. Indeed, Bogomolov's locall Torelli theorem ensures that $\Pi$ is a local diffeomorphism. It is known that $\per$ is simply connected \cite[Claim 2.9]{Ver}. On the other hand, in the next section will see that $\per$ has the path lifting  property for rectifiable paths. This will imply that the period map is a diffeomorphism.

\end{proof}

\subsection{Lifting of generic twistor lines}

\begin{defn}
A twistor sphere $S\in\per$ is called {\it generic} if the corresponding positive $3$-plane $T$ satisfies the following condition: $T^\perp\cap H^2(M,\mathbb{Z})=
\varnothing$. The orthogonal complement is taken here with respect to the Beauville-Bogomolov-Fujiki form.
\end{defn}
The important property of generic twistor spheres is that they can be lifted to the Teichm\"uller space \cite[Proposition 6.1]{Ver}:
\begin{theorem}
Let $S\subset\per$ be a generic twistor sphere and $p\in \teich_b$ a point such that $\Pi(p)\in S$. Then there exists the unique holomorphic sphere $\Sigma\subset \teich_b$, such that $\Pi$ restricted to $\Sigma$ is a diffeomorphism  of $\Sigma$ to $S$.
\end{theorem}

We will call lifts of generic twistor spheres in $\per$ to the Teichm\"uller space also twistor spheres. This allows us to define a metric on $\teich_b$ for which $\Pi$ is a local isometry:

\begin{defn}
Let $p,q\in \teich_b$. Connect them by a chain of twistor spheres $S_i$ such that $S_i\cap S_{i+1}\neq \varnothing$. Compute the inner distance between $p$ and $q$ inside $\cup_i S_i$ with respect to the metric induced from $\per$. Define the distance $d(p,q)$ as an infimum of this number over all chains of twistor spheres connecting $p$ and $q$.
\end{defn}

Note that this metric is subconical. Indeed, in a neighborhood of a point in $\teich_b$ it is isometric to a subtwisor metric on $\per$. Now pick a point $p\in \per$ and take a small ball around it. By Remark \ref{piecewise_spheric} any point of this neighborhood is connected to $p$ by a piecewise smooth path , each smooth part of which lying on some twistor sphere. We can lift this path, thus establishing an isometric section of $\Pi$ over $U$. This proves

\begin{theorem}
The period map $\Pi$ is a covering map.
\end{theorem}

\medskip
\nocite{*}
 
\bibliographystyle{alpha}
\bibliography{biblio.bib}

\end{document}